\documentclass[12pt]{amsart} 
\def\@ptsize{2} 
\newtheorem{Thm}{Theorem}[section]

\newtheorem{Cor}[Thm]{Corollary}
\newtheorem{prop}[Thm]{Proposition}
\theoremstyle{definition}
\newtheorem{rem}[Thm]{Remark}
\begin{document}
\title{On an invariant related to a linear inequality}
\author{Amnon Besser and Pieter Moree}
\address{Department of Mathematics\\Ben-Gurion University of the Negev\\P.O.B. 653\\Beer-Sheva 84105\\Israel}
\address{KdV-Institute\\ University of Amsterdam\\Plantage Muidergracht 24\\ 1018 TV Amsterdam\\ The Netherlands}
\keywords{Linear inequality, invariant, Rademacher function}
\subjclass{Primary 15A39; Secondary 11B99}
\begin{abstract}
Let 
${\underline \alpha}=(\alpha_1,\alpha_2,\ldots,\alpha_m)\in  \mathbb R_{>0}^m$.
Let $\underline \alpha_{i,j}$ be the vector obtained from $\underline \alpha$
on deleting the entries $\alpha_i$ and $\alpha_j$.
We investigate some invariants and near invariants related to the
solutions $\underline \epsilon\in \{\pm 1\}^{m-2}$ of the linear inequality
$|\alpha_i-\alpha_j|<\langle \underline \epsilon,\underline \alpha_{i,j}
\rangle < \alpha_i+\alpha_j$,
where $\langle,\rangle$ denotes the usual inner product. One of
our methods relates, by the use of Rademacher functions, integrals
involving products of trigonometric functions to these quantities.
\end{abstract}
\maketitle\section{Introduction}
The purpose of this note is to construct a certain invariant related to a linear inequality. To give an example, consider the numbers $4,6,7,9$ and $11$. Pick a pair out of them, say $4$ and $6$. Then consider the linear combinations
of the form $\pm 7\pm 9\pm 11$ that are in the open interval
$(|4-6|,4+6)$. There are two of them: $1\cdot 7+(-1)\cdot 9+1\cdot 11$, 
which we give weight $1\cdot (-1)\cdot 1=-1$ and $1\cdot 7+1\cdot9+(-1)\cdot 11$, which we give weight $1\cdot 1\cdot (-1)=-1$. Adding 
yields $-2$. If we now pick any two other numbers and repeat the same construction, we also get $-2$. This is no coincidence; one obtains the same invariance result for any sequence $\alpha_1,\cdots,\alpha_m$ of positive 
reals with $m$ is odd, provided
that $\pm \alpha_1 \pm \alpha_2 \pm \cdots \pm \alpha_m\ne 0$. The value of the invariant depends on the numbers chosen at the outset. If the collection of numbers has even size, we can only obtain an invariance result modulo $2$.\\
\indent A word about the structure of this note: We first give the results and a direct proof that the quantity computed is independent of the choice of the two numbers. Based on our method of proof we can then supply a closed formula for this quantity. At the same time,  having found the closed formula, we can give a much easier alternative proof of the theorem. We 
also give some explanation regarding our choice of weights. Finally,
in Section 3 we rederive some of the results using Rademacher functions
and obtain some new ones.
\section{Results}
\begin{Thm}
\label{main}  
Let $m\ge 3$. Let 
${\underline \alpha}=(\alpha_1,\alpha_2,\ldots,\alpha_m)\in  \mathbb R_{>0}^m$ and suppose that there is no 
$\underline \epsilon \in  \{\pm 1 \}^m$ satisfying $\left< \underline \epsilon,\underline    \alpha \right>=0$. 
Let $1\le i<j\le m$. Let $\underline  \alpha_{i,j}\in 
\mathbb R_{>0}^{m-2}$ be the vector obtained from  $\underline \alpha$ on deleting $\alpha_i$ and $\alpha_j$. Let  
\begin{equation*}    
S_{i,j}(\underline \alpha) :=\{\underline \epsilon\in \{\pm 1    \}^{m-2}: |\alpha_i-\alpha_j|< \quad \left<\underline      \epsilon,\underline \alpha_{i,j}\right> \quad    <\alpha_i+\alpha_j\}.  
\end{equation*}  
\item {\textup{a)}} 
The reduction of $\# S_{i,j}(\underline \alpha)$ mod $2$    only depends on $\underline \alpha$.  
\item {\textup{b)}} Define $N_{i,j}(\underline \alpha) =\sum_{\underline      \epsilon\in S_{i,j}} \prod_{k=1}^{m-2}\epsilon_k$. 
Suppose that    $m$ is odd. Then $N_{i,j}(\underline \alpha)$ only depends on    $\underline \alpha$.
\end{Thm}
\begin{Cor}
Suppose $m=4$. Then $|N_{i,j}(\underline \alpha)|$ only depends
on $\underline \alpha$.
\end{Cor}
\begin{proof} In this case obviously $\# S_{i,j}(\underline \alpha)
\in \{0,1,2\}$. If $\# S_{i,j}(\underline \alpha)\in \{0,2\}$, then $|N_{i,j}(\underline \alpha)|$ is easily seen to equal
zero. If $\# S_{i,j}(\underline \alpha)=1$, then 
$|N_{i,j}(\underline \alpha)|=1$. By property a) it then follows
that $|N_{i,j}(\underline \alpha)|$ is an invariant.
\end{proof}
\begin{rem}
An alternative notation for $N_{i,j}(\underline \alpha)$ and
$S_{i,j}(\underline \alpha)$ we will use that turns out to
be more convenient on occasion is  $N_{\alpha_i,\alpha_j}(\underline \alpha)$, respectively $S_{\alpha_i,\alpha_j}(\underline \alpha)$.
\end{rem}

\noindent {\it Proof of Theorem} \ref{main}.
One checks the result easily in the case $m=3$, so we assume that  $m\ge 4$. Suppose we are given two pairs of indices, $(i_1,j_1)$ and  $(i_2,j_2)$, with one index repeated. We can choose an index $l$ different from all  of 
the above indices. 
We will show that as function of  $\alpha_l$ we have 
$N_{i_1,j_1}(\underline \alpha)=N_{i_2,j_2}(\underline \alpha)$ when the  other components of $\underline \alpha$ are kept fixed, and that the same holds with $N_{\ast,\ast}$ replaced by $\#S_{\ast,\ast}$ modulo $2$. If this holds for all possible pairs with the above restriction, then also $N_{i_1,j_1}(\underline \alpha)=N_{i_2,j_2}(\underline \alpha)$ holds in the case the indices are all different (and the same with $N_{\ast,\ast}$ replaced by 
$\#S_{\ast,\ast}$ modulo $2$).\\ 
\indent Notice that  if $\alpha_l$ is sufficiently large, then $S_{\ast,\ast}(\underline  \alpha) =\{\emptyset\}$. Thus in this case we get an equality with  $\#S_{\ast,\ast}(\underline \alpha)=0$ and $N_{\ast,\ast}(\underline \alpha)=0$.  As $\alpha_l$ decreases, a 
change in some $S_{i,j}(\underline  \alpha)$ will only occur if $\alpha_l$ crosses the point where  $|\alpha_j-\alpha_i|=\left<\underline \epsilon,\underline    \alpha_{i,j}\right>$ or $\alpha_j+\alpha_i= \left<\underline    \epsilon,\underline \alpha_{i,j}\right>$, that is at most at those  $\alpha_l$ such that $\left<\underline \epsilon,\underline    \alpha\right>=0$ for some $\underline \epsilon \in \{\pm 1\}^m$.  We will 
actually see that on moving across such $\underline \alpha$  a change always occurs. In order to prove Theorem \ref{main} it is enough to  prove that both $N_{i,j}(\underline \alpha)$ and  $\#S_{i,j}(\underline \alpha)$ change by the same amount  (respectively the same amount mod $2$), independent of $i,j$, when  we go from  $(\alpha_1,\ldots,\alpha_{l-1},\alpha_l+\delta,\ldots,\alpha_m)$ to  $(\alpha_1,\ldots,\alpha_{l-1},\alpha_l-\delta,\ldots,\alpha_m)$, where  $\left<\underline \epsilon,\underline \alpha\right>=0$ 
and $\delta$  is sufficiently small, but positive.  Without loss of generality we may assume that  $\alpha_j\ge \alpha_i$.  Put $\alpha_l^{+}=\alpha_l+\delta,  \alpha_l^{-}=\alpha_l-\delta$ and, 
for $k\ne l$,  $\alpha_k^{+}=\alpha_k^{-}=\alpha_k$. Let $N_{i,j}^{\pm}(\underline  \alpha)$ and $S_{i,j}^{\pm}(\underline \alpha)$ have the obvious  definitions.  
Note that $\left<\underline \epsilon,\underline    \alpha\right>=0$ implies that  \begin{equation*}    
-\epsilon_j\sum_{k\ne i,j}\epsilon_k\alpha_k^{+}=    \alpha_j^{+}+\epsilon_i\epsilon_j\alpha_i^{+}-\epsilon_j\epsilon_l\delta  \end{equation*}  
and also that  
\begin{equation*}    
-\epsilon_j\sum_{k\ne i,j}\epsilon_k\alpha_k^{-}=    \alpha_j^{-}+\epsilon_i\epsilon_j\alpha_i^{-}+\epsilon_j\epsilon_l\delta.  \end{equation*}  If $\epsilon_i\epsilon_j=1$ and $\epsilon_j\epsilon_l=-1$, the  passage from $\alpha_l^{+}$ to $\alpha_l^{-}$ leads to  $-\epsilon_j\cdot \underline \epsilon_{i,j}$ to be 
added to  $S_{i,j}(\underline \alpha)$. Similarly looking at the other  possible values for $\epsilon_i\epsilon_j$ and  $\epsilon_j\epsilon_l$, we see that for 
all  sign possibilities a  solution is added or deleted according to whether the sign of  $(\epsilon_i\epsilon_j)\cdot  (\epsilon_j\epsilon_l)=\epsilon_i\epsilon_l$  is negative, respectively positive and therefore we conclude  that  $\#S_{i,j}^{-}(\underline \alpha)=\#S_{i,j}^{+}(\underline  \alpha)-\epsilon_i\epsilon_l$.  
In general, let $\underline  \epsilon(1),\ldots,\underline \epsilon(s)$ denote all the different  solutions to  $\left<\underline \epsilon,\underline \alpha\right>=0$,  where $\underline\epsilon$ and $-\underline \epsilon$ are considered  the same solution. Each of 
them leads to a contribution of  $\epsilon_i(r)\epsilon_l(r)$ to $\#S_{i,j}^{-}(\underline  \alpha)-\#S_{i,j}^{+}(\underline \alpha)$ that is not yet  accounted  for, where $1\le r\le s$. If there would be further changes in 
the passage  from $+$ to $-$ they would lead to  additional solutions of $\left<\underline    \epsilon,\underline \alpha\right>=0$. We deduce that $\#  S_{i,j}^{-}(\underline \alpha) =\# S_{i,j}^{+}(\underline \alpha)-  \sum_{r=1}^s \epsilon_i(r)\epsilon_l(r)$. In particular mod $2$ the  latter sum is independent of the choice of $i$ and $j$. This proves  part a.  
Note that $N_{i,j}(\underline \alpha)$ changes by  $\sum_{r=1}^{s} \epsilon_i(r)\epsilon_l(r)\prod_{k\ne    i,j}(-\epsilon_j(r) \epsilon_k(r))$. If $m$ is odd this is equal  to $-\sum_{r=1}^s \epsilon_l(r)\prod_{k=1}^m \epsilon_{k}(r)$. This  is independent of $i$ and $j$, thus proving part b.
\qed\\

\noindent The more interesting part of Theorem \ref{main}, that is part b, raises the question of giving an alternative description of $N_{i,j}(\underline\alpha)$ that does not involve $i$ and $j$. The next theorem will 
give such a description. Let $\operatorname{sgn}(\beta)$ denote the function that equals $1$ if $\beta>0$, $0$ if $\beta=0$ and $-1$ if $\beta<0$.
Let $m\ge 3$ be odd. Let $\delta_0$ be the delta distribution at $0$. This is a generalised function which is the derivative of any step function of jump $1$ at $0$, e.g. $(1/2)\operatorname{sgn}$. When we checked how $N_{i,j}(\underline \alpha)$ varied when we changed $\alpha_l$, we were in fact computing the derivative of $N_{i,j}(\underline \alpha)$ with 
respect to $\alpha_l$. The computation we made in the proof of 
Theorem \ref{main} clearly gives the following formula:
\[  
\frac{\partial N_{i,j}(\underline \alpha)}  {\partial \alpha_l}=-\frac{1}{2} \sum_{\underline  \epsilon \in \{\pm 1\}^m}\epsilon_l\delta_0(\left<\underline  \epsilon,\underline \alpha\right>) \prod_{k=1}^m \epsilon_k.
\]
The factor of $1/2$ comes from the fact that we count each solution together with its negative. This almost proved the next result:
\begin{Thm}
\label{amnon}  
For $m\ge 3$ odd and $\underline \alpha$ as  in Theorem {\rm \ref{main}}, we have  \begin{equation}    
\label{besser}    
N_{i,j}(\underline \alpha)=-{1\over 4}\sum_{\underline \epsilon \in      \{\pm 1\}^m}\operatorname{sgn}(\left<\underline      \epsilon,\underline \alpha\right>)    \prod_{k=1}^m \epsilon_k.  
\end{equation}
\end{Thm}
We give two proofs in this section and another proof, based on Rademacher functions, in Section 3.
\begin{proof}[Proof \textup{1}]  
It follows from the computation above that both sides of  \eqref{besser} have the same partial derivatives. It is also easy to  see that  
\[  \lim_{\alpha_l\rightarrow \infty} \text{R.H.S.} = -{1\over    4}\sum_{\underline \epsilon \in \{\pm 1\}^m}\epsilon_l
\prod_{k=1}^m \epsilon_k =0,  
\]  
and the theorem follows.
\end{proof}
We wish, however, to give a second, more direct proof, which does not use the computation done while proving Theorem \ref{main}.
\begin{proof}[Proof \textup{2}]  
Denote the right hand side in \eqref{besser} by $g(\underline  \alpha)$ and the summand by $h(\underline \alpha,\underline  \epsilon)$. Since $m$ is odd by assumption, we have $h(\underline  \alpha,-\underline \epsilon)= h(\underline \alpha,\underline  \epsilon)$. Thus, we can write  
\[    g(\underline \alpha)=-{1\over 2} \sum_{\left<\underline        \epsilon_{i,j},\underline \alpha_{i,j}\right>>0} h(\underline    \alpha,\underline \epsilon),  
\]  
where we sum over all $\underline  \epsilon\in \{\pm 1\}^m$ satisfying the condition.  (Note that  \[    \sum_{\left<\underline \epsilon_{i,j},\underline \alpha_{i,j}    \right>=0}h(\underline \alpha,\underline \epsilon)=0.)  \]  
From this  
\[    g(\underline \alpha)=-{1\over 2}\sum_{\left<\underline        \epsilon_{i,j},\underline \alpha_{i,j}\right>>0}\left[      \sum_{\epsilon_i,\epsilon_j}\epsilon_i\epsilon_j        \operatorname{sgn}\biggl(\left<\underline        \epsilon_{i,j},\underline        \alpha_{i,j}\right>+\epsilon_i\alpha_i+\epsilon_j\alpha_j\biggr)  \right]    \prod_{k\ne i,j} \epsilon_k.  
\]  
One checks that the  sum in the square brackets is $0$ if $\underline \epsilon_{i,j}$  satisfies either  $\left<\underline \epsilon_{i,j},\underline    \alpha_{i,j}\right>>\alpha_i+\alpha_j$ or 
$\left<\underline    \epsilon_{i,j},\underline \alpha_{i,j}\right><|\alpha_j-\alpha_i|$,  and is $-2$ if $|\alpha_j-\alpha_i|<\left<\underline    \epsilon_{i,j}, \underline\alpha_{i,j}\right><\alpha_i+\alpha_j.$  
Thus  
\begin{equation}
\label{zalmeworstwezen}
g(\underline \alpha)=\sum_{|\alpha_j-\alpha_i|<\left<\underline        \epsilon_{i,j},\underline\alpha_{i,j}\right>      <\alpha_i+\alpha_j}\prod_{k\ne i,j}\epsilon_k=N_{i,j}(\underline    \alpha).  
\end{equation}
(The case where $\left<\underline    \epsilon_{i,j},\underline \alpha_{i,j}\right>  =|\alpha_i-\alpha_j|$ or $\left<\underline \epsilon_{i,j},\underline    \alpha_{i,j}\right> =\alpha_i+\alpha_j$ does not occur, since by  assumption 
$\left<\underline \epsilon,\underline \alpha\right>\ne  0$.)
\end{proof}
\begin{rem}  
Notice that  $g(\underline \alpha)=0$ when $m$ is even.
\end{rem}
Modulo $2$ we have, by (\ref{zalmeworstwezen}), that $g(\underline  \alpha)=\sum_{|\alpha_j-\alpha_i|<\left<\underline      \epsilon_{i,j},\underline\alpha_{i,j}\right>    
<\alpha_i+\alpha_j}1=N_{i,j}(\underline \alpha)$ 
in case $m\ge 3$ is odd.
Thus Theorem \ref{main} a) follows 
in case $m$ is odd. From Theorem \ref{amnon} the validity
of Theorem \ref{main} b) immediately follows.\\
\indent Is there a result similar to Theorem \ref{main} b) when $m$ is even?
Clearly the exact same statement is false. It is however conceivable that there exists an assignment of weights to the elements of 
$S_{i,j}(\underline \alpha)$ that would lead to a similar result. We now show that under certain assumptions 
this is impossible, while also suggesting how one might 
have guessed the correct form of the weights in Theorem \ref{main} in the first place. We consider a weight function $f:\{\pm 1\}^{m-2}\rightarrow \mathbb{R}$. We would like $f$ to be such that $N_{i,j}^\prime(\underline \alpha),$ where
\[   
N_{i,j}^\prime (\underline \alpha):=\sum_{\underline \epsilon\in  
S_{i,j}(\underline \alpha)} f(\underline \epsilon),
\] 
is independent of $i$ and $j$. A look at the proof of Theorem \ref{main} shows that $f$ will have this property if and only if the following condition is satisfied:
\begin{description}
\item[($\ast$)] For any $\underline \epsilon\in \{\pm 1\}^m$ the  quantity $\epsilon_i f(-\epsilon_j\cdot \underline \epsilon_{i,j})$  is independent of $i$ and $j$.
\end{description}
Indeed, since $N_{i,j}^\prime (\underline \alpha)$  is independent of whether $\alpha_j\ge \alpha_i$ or $\alpha_i<\alpha_j,$ we can assume w.l.o.g. that $\alpha_j\ge \alpha_i,$ the proof of Theorem \ref{main} then 
shows that when we pass over a solution of $\left< \underline\epsilon,\underline \alpha\right>=0$,  we gain or lose a solution 
according to the sign of $\epsilon_l\epsilon_i$, with $l$ the coordinate that we vary, and that this solution is $-\epsilon_j\cdot \underline \epsilon_{i,j}$. $N_{i,j}$ therefore changes by $\epsilon_l\epsilon_i f(-\epsilon_j\cdot \underline \epsilon_{i,j})$
and we may neglect $\epsilon_l,$ since it is fixed in the argument. We now have the easy
\begin{prop}  
\label{bestaatniet}
If $m\ge 4$ is even, there is no function $f$ satisfying the condition  \textup{($\ast$)}.  If $m\ge 3$ is odd, every $f$ satisfying \textup{($\ast$)} is of the form  $f(\underline \epsilon)=C\cdot \prod_{k=1}^{m-2} \epsilon_k$, with $C$ a  constant.
\end{prop}
\begin{proof} Let $m\ge 3.$  Consider $(\epsilon_1,\ldots,\epsilon_r,\ldots,\epsilon_{m-2})\in  \{\pm 1\}^{m-2}$. If we apply condition ($\ast$) to the vector  
\[    (\epsilon_1,\ldots,\epsilon_{r-1},\epsilon_r,-\epsilon_r,\epsilon_{r+1},    \ldots,\epsilon_{m-2},-1)\in \{\pm 1\}^m,  
\]  with $(i,j)=(r,m)$ and $(i,j)=(r+1,m),$ we see immediately that $f$  has to satisfy  
\[    f(\epsilon_1,\ldots,\epsilon_r,\ldots,\epsilon_{m-2})=-    f(\epsilon_1,\ldots,-\epsilon_r,\ldots,\epsilon_{m-2}).  
\]  
Therefore, $f(\underline \epsilon)=C\cdot \prod_{k=1}^{m-2}  \epsilon_k$ with $C=f(1,1,\ldots,1)$. Now it is immediately checked  that this function satisfies ($\ast$) only if $m$ is odd.
\end{proof}
\subsection{Shortening vectors}
The quantities above can be related to quantities of the same nature, but
for shortened vectors. Let $\underline \alpha$ be a vector of the
type allowed in Theorem \ref{main}.
For $j\ne k$, let ${\underline \gamma}^{\pm}_{j,k}$ be the vector of
length $m-1$ obtained from 
$\underline \alpha$ on replacing $\alpha_j$ by $|\alpha_j\pm \alpha_k|$
and deleting $\alpha_k$. It can be deduced, for example, that if 
$m\ge 4$ and $\alpha_k\le \alpha_j-\alpha_i$ with $i,j$ and $k$ distinct,
then
\begin{equation}
\# S_{\alpha_i,\alpha_j}(\underline \alpha)=
\# S_{\alpha_i,|\alpha_j-\alpha_k|}({\underline \gamma}^{-}_{j,k})+
\# S_{\alpha_i,\alpha_j+\alpha_k}({\underline \gamma}^{+}_{j,k}).
\end{equation}
To see this, note that if $\alpha_k\le \alpha_j-\alpha_i$
the number of $\underline \epsilon\in
\{\pm 1\}^{m-2}$ with $\langle {\underline \epsilon},
{\underline \alpha}_{i,j}\rangle = \cdots + \alpha_k +\cdots$ and
the number of $\underline \epsilon$ with 
$\langle {\underline \epsilon},
{\underline \alpha}_{i,j}\rangle = \cdots - \alpha_k +\cdots$, equals
$\# S_{\alpha_i,|\alpha_j-\alpha_k|}({\underline \gamma}^{-}_{j,k})$,
respectively 
$\# S_{\alpha_i,\alpha_j+\alpha_k}({\underline \gamma}^{+}_{j,k})$.
We defer further discussion of shortening until Section \ref{kopjekleiner},
where a more powerful approach in uncovering and proving this type
of identities is employed.

\section{Results obtained by using Rademacher functions}
Let $0\le t\le 1$ be a real number. We define $\epsilon_i(t)$, $i=1,2,\cdots$ recursively. Suppose $\epsilon_1(t),\cdots,\epsilon_{i-1}(t)$ are already defined. Then we define $\epsilon_i(t)$ to be $1$ if
$${1\over 2^i}\le t-\sum_{j=0}^{i-1}{\epsilon_j(t)\over 2^j}$$ 
and to be zero otherwise. Note that
$0.\epsilon_1(t)\epsilon_2(t)\cdots$ gives a binary representation for $t$. We define $r_i(t)=1-2\epsilon_i(t)$ for $i=1,2,\cdots$. The functions $r_i(t)$ are called Rademacher functions.\hfil\break
\indent Using the Rademacher functions we will prove the following result,
which shows again that $N_{*,*}(\underline \alpha)$ only depends
on $\underline \alpha$ in case $m$ is odd.
Using (\ref{result}) it is then easy to give yet
another proof of Theorem \ref{amnon}.
\begin{Thm}
\label{pietertje}
Let $m\ge 3$ and $\underline \alpha$ as  be as in Theorem {\rm 1}.
Let $\beta_1,\cdots,\beta_m$ be positive integers and $q$ a real number such
that, for $k=1,\cdots,m$,
$$|{\beta_k\over q}-\alpha_k|<{{\rm min}_{{\underline \epsilon}
\in \{\pm 1\}^m}|\langle {\underline \epsilon},{\underline \alpha}
\rangle|\over m},$$
and with the $\beta$'s satisfying the same ordering and equalities
as do the $\alpha$'s (that is if $\alpha_i\le \alpha_j$, then
$\beta_i\le \beta_j$, where $1\le i,j\le m$ with $i\ne j$).
Suppose $m$ is odd. Then
\begin{equation}
\label{result}
N_{i,j}({\underline \alpha})=(-1)^{m+1\over 2}{2^{m-2}\over 2\pi}\int_0^{2\pi}
\cot({x\over 2})
\sin(\beta_1 x)\cdots \sin(\beta_m x)dx.
\end{equation}
Suppose $m$ is even and $\alpha_i\le \alpha_j$. Then
\begin{equation}
\label{result1}
N_{i,j}({\underline \alpha})=(-1)^{{m\over 2}-1}{2^{m-2}\over 2\pi}\int_0^{2\pi}
\cot({x\over 2})\cot(\beta_j x)
\sin(\beta_1 x)\cdots \sin(\beta_m x)dx.
\end{equation}
Let $m\ge 3$ be arbitrary and $\alpha_i\le \alpha_j$. Then
\begin{equation}
\label{result2}
\# S_{i,j}({\underline \alpha})={2^{m-2}\over 2\pi}\int_0^{2\pi}
\cot({x\over 2})\tan(\beta_i x)
\cos(\beta_1 x)\cdots \cos(\beta_m x)dx.
\end{equation}
\end{Thm}
\begin{Cor}\hfil\break
a) Theorem {\rm \ref{main} b)} holds true.\\
b) If $m\ge 4$ is even, then $N_{*,*}(\underline \alpha)$ depends
only on the largest component omitted from $\underline \alpha$.
This implies that $N_{*,*}(\underline \alpha)$ assumes at most
$m-1$ values.\\
c) The value of $\# S_{*,*}(\underline \alpha)$ depends only
on the smallest component omitted from $\underline \alpha$.\\
\end{Cor}
\noindent {\it Proof of Theorem} \ref{pietertje}.
The existence of $q$ and ${\underline \beta}
=(\beta_1,\cdots,\beta_m)$ is obvious. Clearly the equality
$N_{i,j}({\underline \alpha})=N_{i,j}(q^{-1}\underline \beta)=
N_{i,j}(\underline \beta)$ holds and similarly we have
$S_{i,j}({\underline \alpha})=S_{i,j}({\underline \beta})$.
We only prove the identity (\ref{result}), the proofs of the other
two being very similar.\\
\indent Let $s$ be an integer. Note that
\begin{equation}
\label{flauw}
{1\over 2\pi}\int_0^{2\pi}e^{i s x}dx=\begin{cases}1 &\text{if~}$s=0$;\\0 
&\text{otherwise}.
\end{cases}
\end{equation}
Let us consider 
$N_{m-1,m}(\underline \alpha)=N_{m-1,m}(\underline \beta)$.
Put $m_1=m-2$. We have
$$N_{m-1,m}({\underline \alpha})={2^{m_1}\over 2\pi}\sum_{s=|\beta_{m-1}-\beta_m|+1}^{\beta_{m-1}+\beta_m-1}
\int_0^1 r_1(t)\cdots r_{m_1}(t)\int_0^{2\pi}
e^{i x(\beta_1 r_1(t)+\cdots+\beta_{m_1}r_{m_1}(t)-s)}dx~dt.$$
Every sequence of length $m_1$ of $+1$'s and $-1$'s corresponds to one and 
only on interval $(j/2^{m_1},(j+1)/2^{m_1})$,  with $0\le j\le 2^{m_1}-1$. Thus
$$\int_0^1 r_1(t)r_2(t)\cdots r_{m_1}(t)e^{i x \sum_{k=1}^{m_1}\beta_k r_k(t)}dt=2^{-m_1}\sum_{{\underline \epsilon}\in \{\pm 1\}^{m_1}}(\prod_{k=1}^{m_1}\epsilon_k)e^{ix\sum_{k=1}^{m_1} \epsilon_k\beta_k}$$
$$=i^{m_1}\prod_{k=1}^{m_1}{e^{i\beta_k x}-e^{-i\beta_k x}\over 2i}=i^{m_1}\prod_{k=1}^{m_1}\sin(\beta_k x).$$
Our expression for $N_{m-1,m}({\underline \alpha})$ can thus be rewritten as
$$N_{m-1,m}({\underline \alpha})={(2i)^{m_1}\over 2\pi}\sum_{s=|\beta_{m-1}-\beta_m|+1}^{\beta_{m-1}+\beta_m-1}\int_0^{2\pi} e^{-ixs}\sin(\beta_1 x)\cdots\sin(\beta_{m_1} x)dx.$$
Let us consider the case where $\beta_{m-1}>\beta_m$. Note that $$\sum_{s=\beta_{m-1}-\beta_m+1}^{\beta_{m-1}+\beta_m-1}e^{-ixs}=
e^{-ix\beta_{m-1}}{\sin((\beta_m-{1\over 2})x)\over \sin({x\over 2})}.$$
Since by assumption $m$ is odd, $m_1$ is odd and hence $(2i)^{m_1}$ is purely imaginary. Since, a priori,
$N_{m-1,m}({\underline \alpha})$ is real, we see that we only have to retain the purely imaginary part of $e^{-ix\beta_{m-1}}$, that is $-i\sin(\beta_{m-1}x)$. We thus find that
\begin{equation}
\label{een}
N_{m-1,m}({\underline \alpha})=-i{(2i)^{m_1}\over 2\pi}\int_0^{2\pi}{\sin(\beta_1 x)\cdots \sin(\beta_{m-1}x)\sin(({\beta_m-1/2})x)\over\sin(x/2)}dx.
\end{equation}
Our assumption on $q$ and $\beta_1,\cdots,\beta_m$ implies that
$\langle {\underline \epsilon} , {\underline \beta}\rangle \ne 0$ 
for every ${\underline \epsilon}\in \{\pm 1\}^m$. Thus instead of summing from $s=\beta_{m-1}-\beta_m+1$ to
$s=\beta_{m-1}+\beta_m-1$, we might as well sum from 
$s=\beta_{m-1}-\beta_m$ to $s=\beta_{m-1}+\beta_m$, this then yields \begin{equation} 
\label{twee}
N_{m-1,m}({\underline \alpha})=-i{(2i)^{m_1}\over 2\pi}\int_0^{2\pi}{\sin(\beta_1 x)\cdots \sin(\beta_{m-1}x)\sin(({\beta_m+1/2})x)\over\sin(x/2)}dx.
\end{equation}
On noting that $\sin(\gamma-\delta)+\sin(\gamma+\delta)=2\sin \gamma \cos \delta,$ we finally obtain (\ref{result}) for $i=m-1$ and $j=m$ on adding (\ref{een}) to (\ref{twee}) and averaging. 
By a completely similar reasoning one deals with 
the case where $\beta_{m-1}\le \beta_m$ and one also finds (\ref{result}). Obviously one also arrives at (\ref{result}) if one considers $N_{i,j}(\underline \alpha)$ for arbitrary $1\le i<j\le m$. 
\qed 
\hfil\break

\begin{rem}\hfil\break
a) The condition on the ordering of the $\beta$'s is not needed in
the derivation of (\ref{result}). It is in (\ref{result1}) and
(\ref{result2}) to infer that the assumption $\alpha_i\le \alpha_j$
implies that $\beta_i\le \beta_j$.\\
b) Since (\ref{flauw}) only holds valid for integral $s$, we are
forced to work with the approximation
vector $\underline \beta$, rather than $\underline \alpha$ itself.\\
c) Using that the argument in (\ref{result}) has period $2\pi$ and is
an odd function if $m\ge 2$ and even, one deduces that the integral in (\ref{result}) equals zero in this case.\\ 
\end{rem}

\subsection{The shortening of vectors reconsidered}
\label{kopjekleiner}
The various formulae in Theorem \ref{pietertje}
can be related to each other by invoking
very elementary trigonometric identies such as
$2\sin \alpha \sin \beta =\cos(\alpha+\beta)+\cos(\alpha-\beta)$.
This then yields shortening formulae. In proving them, which is
left to reader, one has to convince oneself that one
can choose an `approximation 
vector' $\underline \beta$ for $\underline \alpha$ that will also 
yield an approximation vector of
the shortened vector(s) involved. An alternative method of proof
is indicated in Section 3.2. 
Recall that ${\underline \gamma}_{j,k}^{\pm}$ is defined in
Section 2.1.
\begin{Thm}
\label{vier}
Let $\underline \alpha$ and ${\underline \alpha}_{i,j}$ be defined
as in Theorem \ref{main}.\\
a) Suppose $m\ge 4$ and even and
$\alpha_i\le \alpha_j$. 
Then
$$N_{i,j}(\underline \alpha)={\rm sgn}(\alpha_j-\alpha_k)N_{*,*}({\underline \gamma}^{-}_{j,k})
-N_{*,*}({\underline \gamma}^{+}_{j,k}).$$
b) Let $m\ge 4$ and $\alpha_i\le \alpha_j$. Suppose furthermore there
exist $r,s\ne i$ such that $\alpha_r+\alpha_s\ge \alpha_i$ and
$|\alpha_r-\alpha_s|\ge \alpha_i$. Then
$$\# S_{\alpha_i,\alpha_j}(\underline \alpha)=
\# S_{\alpha_i,|\alpha_r-\alpha_s|}({\underline \gamma}_{r,s}^{-})+
\# S_{\alpha_i,\alpha_r+\alpha_s}({\underline \gamma}_{r,s}^{+}).$$
c) Let $m\ge 5$ be odd. Suppose we have
$\alpha_k\le |\alpha_i-\alpha_j|$ for some $i,j,k$ with $k\ne i,j$. Then
$$N_{*,*}(\underline \alpha)=N_{\alpha_k,\alpha_i+\alpha_j}
({\underline \gamma}_{i,j}^-)-N_{\alpha_k,|\alpha_i-\alpha_j|}
({\underline \gamma}_{i,j}^{+}).$$
In case $\alpha_i=\alpha_j$ for some $i\ne j$ and $\alpha_k\le 2\alpha_i$ for
some $k\ne i,j$, then
$$N_{*,*}(\underline \alpha)=
-N_{\alpha_k,\alpha_i+\alpha_j}({\underline \gamma}_{i,j}^{+})
-2N_{*,*}({\underline \alpha}_{i,j}).$$
\end{Thm}

\subsection{Theorem \ref{amnon} reconsidered and some analoga}
In this subsection we present a third proof of Theorem \ref{amnon} and
present another theorem that can be proved using the same method of
proof.\\
\noindent {\it Third proof of Theorem} \ref{amnon}.
On inverting some of the last steps in the proof of Theorem \ref{pietertje}, one easily checks that for integer $\beta$ we have
\begin{equation}
\label{drie}
{1\over 2\pi}\int_0^{2\pi}\cot({x\over 2})\sin(\beta x)dx ={\rm sgn}(\beta).
\end{equation}
On writing $\sin(\beta_j x)$ as $(e^{i\beta_j x}-e^{-i\beta_j x})/2i$
for $j=1,\cdots,m$ in (\ref{result}) and
multiplying all these factors out, one
gets a sum of terms of the form
$e^{i\langle {\underline \epsilon},{\underline \beta}\rangle x}\cos(x/2)/\sin(x/2)$, where the term
$e^{i\langle {-\underline \epsilon},{\underline \beta}\rangle x}\cos(x/2)/\sin(x/2)$ appears
with opposite sign, due to the fact that $m$ is odd. This allows one to rewrite (\ref{result}) in the form
$$N_{i,j}(\underline \alpha)=-{1\over 4}\sum_{{\underline \epsilon}\in \{\pm 1\}^m}(\prod_{k=1}^m \epsilon_k){1\over 2\pi}\int_0^{2\pi}\cot({x\over 2})\sin(\langle {\underline \epsilon},{\underline \beta}\rangle x) dx.$$
By (\ref{drie}) we then find the expression
for $N_{i,j}(\underline \alpha)$ as given in Theorem \ref{amnon} with
$\langle {\underline \epsilon},{\underline \alpha}\rangle$ replaced by
$\langle {\underline \epsilon},{\underline \beta}\rangle$.
The proof is now completed on noting that $\underline
\beta$ has the property that
${\rm sgn}(\langle {\underline \epsilon},{\underline \beta}\rangle)
={\rm sgn}(\langle {\underline \epsilon},{\underline \alpha}\rangle)$.
\qed\\

\noindent The latter method of proof can also be applied to equalities (\ref{result1}) and (\ref{result2})
and then yields Theorem \ref{vijf}. Theorem \ref{vijf} is also easily derived on
employing the method of proof in Proof 2 of Theorem \ref{amnon}.
\begin{Thm}
\label{vijf}
Let $m\ge 3$ and $\underline \alpha$ be as in {\rm Theorem \ref{main}}. 
Suppose $\alpha_i\le \alpha_j$.
Then
$$
\# S_{i,j}(\underline \alpha)=
{1\over 2}\sum_{\underline \epsilon\in \{\pm 1\}^m\atop \epsilon_i=1}
{\rm sgn}(\langle \underline \epsilon, \underline \alpha\rangle ).$$
If $m$ is even, we have, moreover,
$$N_{i,j}(\underline \alpha)={1\over 4}
\sum_{\underline \epsilon\in \{\pm 1\}^m}\epsilon_j
{\rm sgn}(\langle \underline \epsilon,\underline \alpha\rangle)
\prod_{k=1}^m \epsilon_k.$$
\end{Thm}
Applying the method of proof of Theorem \ref{amnon} on the right hand side
of the latter identity one finds the following invariant in case $m$ is
even. Note that since $i,j$ are required to be distinct from some prescribed
number, the result is consistent with Proposition \ref{bestaatniet}.
\begin{Thm}
Let $\underline \alpha$ and ${\underline \alpha}_{i,j}$ be as in 
Theorem {\rm \ref{main}}.
Let $h\ne i,j$. Let $h_1$ be an index such that the $h_1$th component
of ${\underline \alpha}_{i,j}$ equlas $\alpha_h$. Define
$$N_{i,j}^{(h)}(\underline \alpha)=
\sum_{{\underline \epsilon}\in S_{i,j}}
\epsilon_{h_1}\prod_{k=1}^{m-2}\epsilon_k.$$
Then, for $m\ge 4$ with $m$ even, $N_{i,j}^{(h)}(\underline \alpha)$ does not depend
on $i$ and $j$. If $\alpha_{h_2}\le \alpha_h$ for some $h_2\ne h$,
then $N_{i,j}^{(h)}(\underline \alpha)=N_{h_2,h}(\underline \alpha)$.
\end{Thm} 
On invoking Theorems \ref{amnon} and \ref{vijf} it is not difficult to reprove Theorem \ref{vier}.
This is left to the reader.
\subsection{An example}
Let $p_1,p_2,\cdots$ denote the consecutive primes. Since for
the natural integers we have
unique factorisation (up to order of factors), we have that $\pm \log p_1 \pm \log p_2 \pm
\cdots \pm \log p_n \ne 0$ and hence we can apply Theorem \ref{main}
with ${\underline \alpha}^{(n)}=(\log 2,\cdots,\log p_n)$.
It is
not difficult to show that, for $1\le i<j\le n$,
$$\# N_{i,j}(\underline \alpha^{(n)})=
(-1)^n \sum_{\sqrt{p_1\cdots p_n}/p_i < m < \sqrt{p_1\cdots p_n}
\atop (m,p_ip_j)=1,~P(m)\le p_n}\mu(m),$$
where $P(m)$ denotes the largest prime factor of $m$ and $\mu$ the
M\"obius function.
If $n\ge 3$ is odd, then the latter quantity does not depend on
$i$ and $j$ by Theorem \ref{main} b).
If $n\ge 4$ is even, the latter quantity does not depend
on $i$ by Corollary 3.2 b). The values one finds of
$N_{*,*}(\alpha^{(n)})$ for $n=5,7,11,13,15$ are, respectively,
$-1,3,22,-53,-55$.\\
\hfil\break
\noindent {\bf Acknowledgement}. Alexander Reznikov suggested during his stay at the Max-Planck-Institute some time in 1997 that preliminary computations 
(joint with Luca Migliorini) in the cobordism theory of the moduli space of polygons suggested that there should be an invariant lurking around, the simplest candidate being 
$\# S_{i,j}(\underline \alpha)$. Our intuition was that the latter clearly could not be an invariant and that the existence of 
a more complicated invariant was highly unlikely. In setting out to establish this, we found to our surprise that there is indeed an associated 
invariant ($N_{i,j}(\underline \alpha)$), though alas $\# S_{i,j}(\underline \alpha)$ does not quite qualify for this r\^ole...\\
\indent We like to thank A. Reznikov for suggesting the problem and the Max-Planck-Institute for its support and scientifically rewarding atmosphere.
\end{document}